\documentclass{amsart}
\usepackage{amsmath}%
\usepackage{amsfonts}%
\usepackage{amssymb}%
\usepackage{graphicx}
\usepackage{mathtools, tikz, float,hyperref,authblk}
\usepackage{amsaddr}

\usepackage{color}

\newcommand{\B}{\mathcal{B}}
\newcommand{\C}{\mathcal{C}}

\newcommand{\PC}{\mathbb{P}\mathcal{C}}

\newcommand{\G}{\mathcal{G}}
\newcommand{\h}{\mathbb{H}^2}
\newcommand{\M}{\mathcal{M}}

\newcommand{\R}{\mathbb{R}}

\newcommand{\T}{\mathcal{T}}

\newcommand{\PML}{\mathcal{PML}}
\newcommand{\ML}{\mathcal{ML}}

\DeclareMathOperator{\Mod}{Mod}

\DeclareMathOperator{\sys}{sys}
\DeclareMathOperator{\Col}{Col}

\newtheorem{theorem}{Theorem}[section]
\theoremstyle{plain}
\newtheorem{lem}[theorem]{Lemma}

\newtheorem{cor}[theorem]{Corollary}
\newtheorem{rem}[theorem]{Remark}
\newtheorem{prop}[theorem]{Proposition}

\newtheorem{ques}{Question}

\theoremstyle{definition}
\newtheorem{defi}[theorem]{Definition}

\long\def\symbolfootnote[#1]#2{\begingroup%
\def\thefootnote{\fnsymbol{footnote}}\footnote[#1]{#2}\endgroup}

\title{A projection from filling currents to Teichm\"uller space}
\author[S. Hensel]{Sebastian Hensel$^{\dag}$}
\address{Mathematisches Institut der LMU\\
	Theresienstr. 39\\
	D-80333 M\"unchen, Germany}
\email{hensel@math.lmu.de}
\thanks{$\dag$ Supported by DFG SPP 2026 “Geometry at Infinity”}

\author[J. Sapir]{Jenya Sapir}
\address{Department of Mathematics and Statistics\\
Binghamton University\\
4400 Vestal Parkway E\\
Binghamton, New York 13902,
USA}
\email{sapir@math.binghamton.edu}

\begin{document}

\begin{abstract}
 Let $S$ be a closed, genus $g$ surface. The space of geodesic currents on $S$ encompasses the set of closed curves up to homotopy, as well as Teichmuller space, and many other spaces of structures on $S$. We show that one can define a mapping class group equivariant, length-minimizing projection from the set of filling geodesic currents down to Teichmuller space, and prove some basic properties of this projection to show that it is well-behaved.
\end{abstract}

\maketitle

\section{Introduction}
Let $S$ be a closed genus $g$ surface and let $\T(S)$ denote its
Teichm\"uller space.  A common theme in Teichm\"uller theory is to try
to minimise the length of a (sufficiently complicated) object on $S$
over all of $\T(S)$, and to look for a metric at which the length is minimized. This idea underlies Kerckhoff's proof of the
Nielsen realisation problem \cite{Kerckhoff83} (where the orbit of a large enough
collection of curves is minimised), and his \emph{lines of minima} \cite{Kerckhoff92},
which yield quasigeodesics in $\T(S)$ with interesting properties \cite{CRS06,CRS07} (where the sum of the lengths of two laminations that together fill the surface is minimised).
Minimising the length (and estimating this minimum) for non-simple closed curves on $S$
has also recently attracted some attention both in its own right \cite{Yamada,Nakanishi,Basmajian13,AGPS16,Gaster17} and as a way to count mapping class group orbits of such curves with fixed self-intersection number \cite{SapirOrbits,AS16}. 

\medskip The purpose of this paper is to define and begin to study a framework generalising
this minimisation idea in the realm of \emph{geodesic currents}, which are in
themselves a powerful tool for studying $\T(S)$. 

The space of geodesic currents $\C(S)$ on $S$ (defined carefully
below) can be viewed as the space of geodesic flow invariant measures
on the unit tangent bundle $T_1(S)$, where we can choose any
hyperbolic metric on $S$ to make this precise. The space of geodesic currents was first defined by Bonahon, and shown to contain $\T(S)$ as well as the set of homotopy classes of closed curves on $S$ \cite{Bon88}. It was later shown to  contain many other geometric
structures such as the space of all negatively curved metrics \cite{Otal90} and the space of flat structures on $S$ \cite{DLR10}. All of these structures are
unified by an intersection function $i(\cdot, \cdot)$ on pairs of
currents. When two currents $\mu$ and $\nu$ represent closed curves on
$S$, then $i(\mu, \nu)$ is just the geometric intersection
number. That is, it is just the minimal number of intersections
between elements in the free homotopy classes of the two
curves. However, if $\mu$ is still a closed curve and $\nu$ represents
a metric, then $i(\mu, \nu)$ becomes the geodesic
length of $\mu$ with respect to the metric represented by $\nu$. 

As a fairly immediate consequence of work of Wolpert
\cite{Wolpert06}, we prove in Proposition \ref{Minimum_exists} that one can define a mapping class group equivariant 
projection
\[
\pi : \C_{fill}(S) \to \T(S)
\]
from the subset of ``filling'' currents $\C_{fill}(S)$ (precisely
defined below) to Teichmuller space.  The map $\pi$ is uniquely
determined by the property that if $\mu \in \C_{fill}(S)$ and
$\pi(\mu) = X_\mu$, then for all $Y \in \T(S)$,
\[
 i(\mu, X_\mu) \leq i(\mu, Y).
\]

 In particular $\pi$ generalises the length-minimisation problem (for
 $\mu$ the counting current of a filling curve, or the sum of currents
 of a filling collection of curves), and the definition underlying
lines of minima (for $\mu$ the sum of two filling measured
laminations).

\medskip
The space of currents has a natural action by $\R_+$, given by scaling each measure by a positive constant. Taking the quotient of $\C(S)$ by this action gives the space $\PC(S)$ of projective currents. The intersection form is bilinear, and so our projection map also yields
\[
 \pi: \PC_{fill}(S) \to \T(S)
\]
where $\PC_{fill}(S)$ is the set of projectivized \textit{filling} currents.
Our main result is the following: 
\begin{theorem}
 The map $\pi: \PC_{fill}(S) \to \T(S)$ is mapping class group equivariant, continuous and proper.
\end{theorem}

In \cite{Kerckhoff83, Kerckhoff92} Kerckhoff shows that when the sum of a collection of simple closed curves or measured laminations gives a filling geodesic current, then the sum has a unique length minimizing metric. Moreover, he shows that, on certain 6g-6-dimensional subspaces of $\PC_{fill}(S)$, the length minimizing projection is mapping class group equivariant, continuous and proper. Although not explicitly written down, the proofs in \cite{Kerckhoff92} can be used to prove Theorem 1.1. We give an independent proof of this result, using more geometric techniques, and relying on properties of currents.

As was pointed out to the authors by Marc Burger, this result is connected to a conjecture in higher Teichm\"uller theory. Higher Teichmuller spaces are spaces of representations from $\pi_1(S)$ to higher rank simple Lie groups. An important problem in higher Teichmuller theory is to find a mapping class group invariant complex structure on these spaces. Labourie conjectured that there is a certain $\Mod(S)$-equivariant projection from higher Teichmuller spaces to $\T(S)$ that induces such a structure \cite[Conjecture 1.0.8]{Labourie}. As shown, for example, in \cite{MZ, BIPP2}, there are many higher Teichmuller spaces that can be embedded into $\C(S)$ as filling currents.  As a consequence, these higher Teichm\"uller spaces then equivariantly surject onto (usual) Teich\"uller space via $\pi$. 
However, the relation to the surjection conjectured by Labourie is unclear, as is the question if it can be used to obtain a complex structure on these spaces.

\smallskip
By work of Diaz and Series, explained in more detail in Section \ref{sec:Extention_to_non-filling}, the projection $\pi$ does not extend to the boundary. However, we extend a result of theirs to show that the projection does extend at maximal, uniquely ergodic laminations in a surprising way. In fact, $\T(S)$ is embedded in $\PC(S)$ in such a way that its boundary is exactly the Thurston boundary. Moreover, $\PML(S)$ also lies in the boundary of $\PC_{fill}(S)$.  We show that the projection extends to maximal, uniquely ergodic laminations, and that, in fact, it extends as the identity map at those points:
\begin{prop}
\label{prop:Extension_at_Uniquely_Ergodic}
Suppose $\lambda$ is a maximal, uniquely ergodic lamination. If $\mu_n \in \C_{fill}(S)$ with $\mu_n \to \lambda$, then $\pi(\mu_n) \to [\lambda] \in \PC(S)$.
Conversely, if $\mu_n \in \C_{fill}(S)$ with $\mu_n \to \mu$ for some $\mu$, and $\pi(\mu_n) \to \lambda$, then, in fact, $[\mu] = [\lambda] \in \PC(S)$.
\end{prop}
\subsection{Background}
We now recall the definition of the space of geodesic currents $\C(S)$ on $S$. Given $X \in \T(S)$, we can identify the universal cover of $X$ with $\h$, and write $X = \h/ \Gamma$ where $\Gamma$ is the image of a representation of $\pi_1(S)$ into $PSL(2,\R)$. The space of geodesic currents is the space of locally finite, $\Gamma$-invariant Borel measures on the set of geodesics in $\h$. We endow $\C(S)$ with the weak-* topology. Note that, a priori, $\C(S)$ depends on the choice of metric $X$. However, for any two metrics $X$ and $Y$, there is a Holder continuous homeomorphism between the sets of currents defined with respect to these metrics.

Note that $\C(S)$ is closed under multiplication by positive real scalars, and under addition of currents. Quotienting $\C(S) \setminus \{0\}$ by the $\R_+$ action gives us the set $\PC(S)$ of projective geodesic currents. Bonahon shows that $\PC(S)$ is compact in \cite{Bon88}.

The set of closed geodesics on $X$ embeds into $\C(S)$ in the following way: given any geodesic $\gamma$ on $X$, we can take its full preimage, $\tilde \gamma$, in $\h$. If $\gamma$ is closed, then $\tilde \gamma$ is a discrete set. Thus, the Dirac measure on $\gamma$ is a geodesic current. One can similarly embed the set $\ML$ of measured laminations into $\C(S)$. Abusing notation, if $\gamma$ is a closed geodesic or measured lamination, we will also use $\gamma$ to describe the corresponding current. The set of weighted closed (not necessarily simple!) geodesics is dense in $\C(S)$ \cite{Bonahon85}. 

Teichmuller space also embeds into the space of geodesic currents. Given $X$, the Liouville measure on $\h$ defines a $\Gamma$-invariant measure on the set of geodesics in $\h$. This is the current corresponding to $X$. Any other $Y \in \T(S)$ defines a homeomorphism $\phi: S^1 \to S^1$ on the boundary at infinity of $\h$ by lifting and extending the homeomorphism from $X$ to $Y$ coming from the markings. Pulling the Liouville current back along $\phi$ gives the geodesic current corresponding to $Y$. Bonahon showed that the resulting map from $\T(S)$ to $\C(S)$, and also the projection to $\PC(S)$ are homeomorphisms onto their image. Again abusing notation, if $X \in \T(S)$, then we will also use $X$ to denote the Liouville current corresponding to $X$.

Lastly, Bonahon showed that the geometric intersection number on geodesics extends continuously to a symmetric, bilinear function on $\C(S) \times \C(S)$. Moreover, if $X \in \T(S)$ and $\gamma \in \G$, then 
\[
 i(\gamma, X) = \ell_X(\gamma)
\]
where $\ell_X(\gamma)$ is the length of the geodesic representative of $\gamma$ on $X$.

We define the \textbf{systolic length} function by
\[
 \sys : \C(S) \to  \R_{\geq 0}, 
\]
\[
\sys(\mu) =  \inf_\alpha i(\mu, \alpha)
\]
where the infimum is taken over all simple closed curves $\alpha$. Let $\mu \in \C(S)$ be a geodesic current. 

We say that a current $\mu$ is \textbf{filling} if $i(\mu, \lambda) > 0$ for all $\lambda \in \C(S)$. Equivalently, every geodesic in $\h$ intersects some geodesic in the support of $\mu$ transversely. Note that if $\mu$ is a (non-simple) closed geodesic, then $\mu$ is filling if and only if it cuts $S$ into simply connected regions. In \cite[Theorem 4.1]{BIPP}, they show that a current $\mu$ is filling if and only if $\sys(\mu) > 0$. This is the characterization of filling we use in this paper. 

We also recall that Bonahon \cite{Bonahon86} has shown that a current $\mu$ is a measured lamination if and only if $i(\mu, \mu) = 0$ (in particular, measured laminations are never filling). 

Observe that a current $\mu$ is filling if and only if $c\mu$ is filling for any $c>0$. We let $\C_{fill}(S)$ denote the set of filling currents, and $\PC_{fill}(S)$ denote the image of this set in $\PC(S)$. Note that a current is filling if and only if its image lies in $\PC_{fill}(S)$.

\section{Projection onto Teichmuller space}
We will show that there is a well-defined projection from filling geodesic currents to Teichmuller space.

\begin{prop}
\label{Minimum_exists}
 If $\mu$ is a filling current, there is a unique point $X \in \T(S)$ so that for all $Y \neq X \in \T(S)$,
 \[
  \ell_X(\mu) < \ell_Y(\mu)
 \]

\end{prop}
\begin{proof}

 For each $\mu \in \C_{fill}(S)$, we define
\[
 \ell_{min}(\mu) =  \inf_{Y \in \T(S)} \ell_Y(\mu)
\]
We will show that there is a unique point $X \in \T(S)$ which realizes $\ell_{min}(\mu)$.

Fix a point $X_0 \in \T(S)$. Let $L = \ell_{X_0}(\mu)$. Then, by definition, $\ell_{\min}(\mu) \leq L$. Consider the set 
\[
 M = \{Y \in \T(S) \ | \ \ell_Y(\mu) \leq L\}
\]
If $\mu$ has a length minimizer, then it must be in this set. 

Similarly, consider the set 
\[
  M' = \{ \nu \in \C(S) \ | \ i(\nu, \mu) \leq L\}
\]
of currents with bounded intersection with $\mu$. By
\cite[Proposition~4]{Bon88}, the set $M'$ is compact in the space
$\C(S)$ of currents on $S$. Further recall that we can properly embed
$\T(S)$ into $\C(S)$ by assigning to a point in Teichm\"uller space
the Liouville current of the corresponding metric \cite[Lemma 10,
Corollary 11]{Bon88}. By construction, we have $M = M' \cap \T(S)$
with respect to this embedding, and thus $M$ is compact. Hence, $\mu$
has a length-minimiser in $M$.

But a minimum in $M$ must be a global minimum on all of $\T(S)$, as
$M$ is a sub-level set of the length function of $\mu$. By Wolpert
\cite{Wolpert06}, the hyperbolic length of a current is a strictly convex function with
respect to the Weil-Petersson metric. Thus, the metric on which the length of $\mu$ is minimised is unique.
\end{proof}

For any current $\mu$ and points $X,Y \in \T(S)$, we have $i(X, \mu) < i(Y, \mu)$ if and only if $i(X, t\mu) < i(Y, t\mu)$ for any $t > 0$. Thus, this projection is actually defined on $\PC_{fill}(S)$:

\begin{defi}
 We define a map 
 \[
  \pi : \PC_{fill}(S) \to \T(S)
 \]
where, for each $[\mu] \in \PC_{fill}(S)$, $\pi[\mu] = X$ if and only if $i(X, \mu) < i(Y, \mu)$ for all $Y \neq X \in \T(S)$. By the remark above, this definition is independent of the choice of representative $\mu$ of $[\mu]$.
\end{defi}
Note that the projection is also mapping class group equivariant, since the action of the mapping class group on $\C(S)$ preserves intersection number.

\begin{rem}
	Bonahon shows in \cite[in Theorem 19]{Bon88} that $i(X,X) \leq i(X,Y)$ for all $X,Y \in \T(S)$. Since the length function is strictly convex, this shows that in fact, $i(X,X) < i(X,Y)$ for all $X \neq Y \in \T(S)$. Hence, the map 
	$\pi$ is constant on $\T(S)$, and therefore actually is a projection. 
\end{rem}

We note that it is also possible to prove the existence of a length minimizer for filling currents using the following extension of the collar lemma, which we include for completeness.

\begin{lem}
\label{lem:CollarLemma}
 Suppose $\mu \in \C(S)$ and $X \in \T(S)$. Let $\alpha$ be a simple closed curve. Then,
 \[
  i(X, \mu) \geq \Col_X(\alpha) i(\mu, \alpha)
 \]
 where $\Col_X(\alpha) = \sinh^{-1}(1/\sinh(\frac 12\ell_X(\alpha))$ is the width of an embedded collar neighborhood around $\alpha$.
\end{lem}
\begin{proof}
 First, take $\gamma$ to be a closed geodesic. Then by the collar lemma, 
 \[
  \ell_X(\gamma) \geq \Col_X(\alpha) i(\gamma, \alpha)
 \]
 Since $\ell_X(\gamma) = i(X, \gamma)$, and the intersection form is linear, the same equation is true if we replace $\gamma$ by $c \gamma$ for any $c > 0$.
 
Now take any $\mu \in \C(S)$. Since weighted closed geodesics are dense in $\C(S)$, there is a sequence $c_n \gamma_n$ where $c_n > 0$ and $\gamma_n$ is a closed geodesic, so that $\lim c_n \gamma_n = \mu$. Then for all $n$ we have 
\[
 i(X,c_n \gamma_n) \geq \Col_X(\alpha) i(c_n \gamma_n, \alpha)
\]
Since the intersection function is continuous, we can take limits of both sides and get the desired inequality.
\end{proof}

\section{Continuity of the projection map}
In this section, we will show that the projection map is continuous. First, we need the following lemma:
\begin{lem}
\label{lem:Minimum_Length}
 Let $\mu$ be a filling current.  Then we have
 \[
  \ell_{\pi(\mu)}(\mu) \leq 4\sqrt{2 |\chi(S)|} \sqrt{i(\mu, \mu)}
 \]
\end{lem}
 For closed geodesics, this result was first proven by \cite{AGPS16} for some constant, and then \cite{AS16} showed that we can take the constant to be $4\sqrt{2 |\chi(S)|}$. We just need to show that we can extend this result to all filling currents.
\begin{proof}
 Let $\mu$ be a filling current. Since closed geodesics are dense in $\PC(S)$, there is a sequence of closed geodesics $\gamma_n$ and constants $c_n > 0$ so that 
 \[
  c_n \gamma_n \to \mu
 \]
Since $\mu$ is filling, $\sys(\mu) > 0$. Since the systole function is continuous, $\sys(c_n \gamma_n) > 0$ for all $n$ big enough. In particular, up to passing to a subsequence, we can assume that $\gamma_n$ are all filling closed curves.

Let $X_n = \pi(\gamma_n)$. This also means that $X_n = \pi(c_n \gamma_n)$. The compactification of $\T(S)$ in $\PC(S)$ is just the Thurston compactification (by Bonahon \cite{Bon88}.) Thus, either we can pass to a subsequence so that $X_n \to Y$ for some $Y \in \T(S)$, or there are constants $d_n > 0$ so that $d_n X_n \to \lambda$ for $\lambda \in \ML(S)$.

Suppose first that $d_n X_n \to \lambda \in \ML(S)$. Then on the one hand,
\[
 \lim_{n \to \infty} i(c_n \gamma_n, d_n X_n) \to i(\mu, \lambda)
\]
Since $\mu$ is filling, by Theorem 4.1 of \cite{BIPP}, we must have $i(\mu, \lambda) > 0$.

However, $i(c_n \gamma_n, c_n \gamma_n) \to i(\mu, \mu)$. Since $i(\mu, \mu) > 0$ and $i(\gamma_n, \gamma_n) \neq 0$, we can replace each $c_n$ with $c_n = \frac{\sqrt{i(\mu, \mu)}}{\sqrt{i(\gamma_n, \gamma_n)}}$ so that convergence still holds. 
So,
\begin{align*}
 \ell_{X_n}(c_n \gamma_n) & \leq c_n 4\sqrt{2 |\chi(S)|} \sqrt{i(\gamma_n, \gamma_n)} \\
  & = 4\sqrt{2 |\chi(S)|} \sqrt{i(\mu, \mu)}
\end{align*}
In particular, $\ell_{X_n}(c_n \gamma_n)$ is bounded. Moreover, we know that $i(X_n, X_n) = 2 \pi |\chi(S)|$ and $i(\lambda, \lambda) = 0$ by \cite{Bon88}. Thus, we must have $d_n \to 0$. So, 
\[
  \lim_{n \to \infty} i(c_n \gamma_n, d_n X_n) = 0
\]
This contradicts the above conclusion. 

Therefore, up to passing to a subsequence, there is some $Y \in \T(S)$ so that $X_n \to Y$. Thus, since 
\[
 i(X_n, c_n \gamma_n) \leq 4\sqrt{2 |\chi(S)|} \sqrt{i(\mu, \mu)}
\]
we can take the limit of the left hand side as $n$ goes to infinity to get that 
\[
 i(Y, \mu) \leq 4\sqrt{2 |\chi(S)|} \sqrt{i(\mu, \mu)}
\]
Since $\pi(\mu)$ is the length minimizer of $\mu$, we have $i(\pi(\mu),\mu) \leq i(Y, \mu)$, giving us the needed result. (In fact, we prove below that $Y = \pi(\mu)$, but we don't need that here.)
\end{proof}

We remark that in \cite{AS16} it is shown that for every non-simple closed curve $\gamma$, there is some metric $X$ so that $\ell_X(\gamma) \leq 4 \sqrt{2 |\chi(S)|} \sqrt{i(\gamma, \gamma)}$. Although the proof above crucially
depends on the fact that $\mu$ is filling, this motivates the following:
\begin{ques}
	Does Lemma~\ref{lem:Minimum_Length} hold for all currents with non-zero self-intersection?
\end{ques}
We are now ready to prove the first part of the main theorem:
\begin{lem}
 The projection $\pi : \PC_{fill}(S) \to \T(S)$ is continuous.
\end{lem}
\begin{proof}
Suppose $\{\mu_n\}$ is a sequence of filling currents that converges to a filling current $\mu$. Let $X_n = \pi(\mu_n)$ be the length minimizer of $\mu_n$ and let $X = \pi(\mu)$ be the length minimizer of $\mu$. First, we can run essentially the same argument as above to show that, up to passing to a subsequence, $X_n \to Y$ for some $Y \in \T(S)$.

Indeed, suppose there are constants $d_n > 0$ so that $d_n X_n \to \lambda \in \ML(S)$. Then 
\[
 \lim_{n \to \infty} i(\mu_n, d_n X_n) = i(\mu, \lambda) > 0
\]
since $\mu$ is a filling current. 

However, $\ell_{X_n}(\mu_n) \leq 4\sqrt{2 |\chi(S)|} \sqrt{i(\mu_n, \mu_n)}$ for each $n$. Since $i(\mu_n, \mu_n) \to i(\mu, \mu)$, the sequence $i(\mu_n, X_n)$ is bounded. Moreover, $\lambda \in \ML(S)$ implies $d_n \to 0$. Thus, 
\[
 \lim_{n \to \infty} i(\mu_n, d_n X_n) = 0
\]
which gives the same contradiction as above.

Thus, after passing to a subsequence, $X_n \to Y$ for some $Y \in \T(S)$. But 
\[
 i(X_n, \mu_n) < i(X, \mu_n)
\]
for each $n$, since $X_n$ is the length minimizer of $\mu_n$. Taking limits, we see that 
\[
 i(Y, \mu) \leq i(X, \mu)
\]
Since $X$ is the length minimizer of $\mu$, this implies that $Y = X$. Thus, $X_n \to X$.
\end{proof}

\section{Non-filling currents}
Before we continue, we will show that non-filling currents do not have
a minimizer in $\T(S)$, extending a result of Kerckhoff on non-filling sums of
laminations \cite{Kerckhoff92}. We highlight that as a consequence, the definition of
the projection
$\pi$ does not extend to a larger class of currents. For a full discussion, see Section \ref{sec:Extention_to_non-filling} below.
\begin{lem}
\label{lem:No_Minimizers_Nonfilling}
 Let $\mu \in \C(S)$ be non-filling. Then there is no $X \in \T(S)$ so that $\ell_X(\mu) \leq \ell_Y(\mu)$ for all $Y \in \T(S)$.
\end{lem}
\begin{proof}
 Let $\mu$ be non-filling. By \cite{BIPP}, either $\mu$ is a measured lamination, or we can find a non-empty collection of simple closed curves $\delta_1, \dots, \delta_n$ so that $\mu = \sum c_i \delta_i + \sum_F \nu_F$, where $c_i > 0$ and the second sum is over connected components $F$ of $S \setminus \cup \delta_i$, and the support of $\nu_F$ is $F$.

If $\mu$ is a measured lamination, then by \cite[proof of Thm 2.1, part 2]{Kerckhoff92} the function $\ell_\mu : \T(S) \to \R$ sending $X \in \T(S)$ to $\ell_X(\mu)$ has no critical points. That is, the length of $\mu$ is not minimized inside $\T(S)$.

So suppose 
 \[
  \mu = \sum c_i \delta_i + \sum_F \nu_F
 \]
 where the set of simple closed curves $\{\delta_1, \dots, \delta_n\}$ is non-empty, and suppose $\mu$ is minimized at a point $X \in \T(S)$. We will build a metric $Y$ so that $\ell_Y(\mu) \leq \ell_X(\mu)$. 
 
 Take $\delta = \delta_1$, and let $Z$ be the metric we get by cutting $X$ along $\delta$. Let $\delta_+$ and $\delta_-$ be the two boundary components of $Z$ coming from $\delta$.
 
 In \cite{PT09}, they give a construction that, for all $\epsilon > 0$ small enough, produces a new metric $Z'$ so that $\ell_Z(\gamma) \leq \ell_{Z'}(\gamma)$ for all simple closed curves $\gamma$, and the length of $\delta$ decreases by $\epsilon$. Since \cite{PT09} only prove this result for simple closed curves, we will summarize their construction here and show that their results give us the statement we want. (Note that a similar result was first proven by Parlier in \cite{Parlier05}, but the construction in \cite{PT09} is better suited to our purposes.)
 
 Let $\alpha$ be a simple arc from $\delta_+$ to itself that is orthogonal to $\delta_+$ at its endpoints. Let $\hat Z$ be the Nielsen completion of $Z$, formed by attaching infinite cuffs to all boundary curves of $Z$. Extend $\alpha$ to an infinite geodesic arc. Papadopoulos and Th\'eret show how to choose two geodesic arcs hyperparallel to $\alpha$ that cut out a neighborhood $B$ about $\alpha$. By cutting out this neighborhood and gluing the two sides in a prescribed way, they obtain a new complete hyperbolic metric $\hat Z_B$. They show that the induced map 
 \[
  f: \hat Z \to \hat Z_B
 \]
 is 1-Lipschitz. Thus, given any closed geodesic $\gamma$ on $\hat Z$, its image $f(\gamma)$ has length bounded above by $\ell_Z(\gamma)$. Let $Z_B$ be the convex core of $\hat Z_B$. Since the geodesic representatives of all closed geodesics lie in the convex core, we have that $\ell_Z(\gamma) \geq \ell_{Z_B}(\gamma)$. By bringing the curves about $\alpha$ closer together, we see that for all $\epsilon > 0$ small enough, we can get $\ell_{Z_B}(\delta_+) = \ell_Z(\delta_+) - \epsilon$. Moreover, the length of $\delta_-$ remains unchanged.
 
 We can repeat this procedure with $\delta_-$, taking care that it gets shrunk by the same amount as $\delta_+$. This will give us a new metric where the lengths of all closed curves do not increase, and the lengths of $\delta_+$ and $\delta_-$ are equal and strictly smaller than before.
 
 Note that $\mu - c_1 \delta$ can be viewed as a current on $S \setminus \delta$. Moreover, since $S \setminus \delta$ is a compact surface with boundary, work of Bonahon \cite{Bonahon85} implies that weighted closed curves are dense in the space of currents. Thus, the fact that the lengths of all closed curves do not increase implies that the lengths of all currents supported on $S \setminus \delta$ do not increase, as well.
 
 We then glue $\delta_+$ to $\delta_-$ to get a new metric $Y \in \T(S)$. Then $Y$ satisfies 
 \[
  \ell_Y(\mu) < \ell_X(\mu)
 \]
as $\ell_Y(\mu - c_1 \delta_1) \leq \ell_X(\mu - c_1 \delta_1)$, and $\ell_Y(\delta_1) < \ell_X(\delta_1)$. 
This is a contradiction, since we assumed that $X$ is a length minimizer.
\end{proof}

\section{Minimal length bound}

We wish to understand the length of a current at points near its length minimizer. Given $X, Y \in \T(S)$, recall that the (asymmetric) Thurston distance is defined by
\[
 d_{Th}(X,Y) = \log \inf \{L \ | \ f: X \to Y \text{ is } L\text{-Lipschitz}\}
\]
Then we have 
\begin{lem}
\label{lem:Lipschitz_length}
 Let $X, Y \in \T(S)$ and let $\mu \in \C(S)$. Then, 
 \[
  i(Y, \mu) \leq e^{d_{Th}(X,Y)} i(X,\mu)
 \]
\end{lem}
\begin{proof}
 Let $X, Y \in \T(S)$. Thurston showed \cite{Thurston88} that the distance from $X$ to $Y$ is realized by a Lipschitz map. So let $f: X \to Y$ be an $L$-Lipschitz map so that 
 \[
  d_{Th}(X,Y) = \log L
 \]
If $\gamma$ is a closed geodesic with respect to $X$, then $f(\gamma)$ is a closed curve of length at most $L \ell_X(\gamma)$. Therefore, 
\[
 i(Y, \gamma) \leq L i(X,\gamma)
\]
where $L = e^{d_{Th}(X,Y)}$.

Now let $\mu \in \C(S)$. We can approximate $\mu$ by a sequence $c_n \gamma_n$ with $c_n > 0$ and $\gamma_n$ closed curves. Then the desired result follows from linearity and continuity of intersection number.
\end{proof}

\section{Properness}
We are now ready to show that the second part of the main theorem.

\begin{theorem}
 The projection $\pi : \PC_{fill}(S) \to \T(S)$ is proper.
\end{theorem}
\begin{proof}
 Since $\T(S)$ is locally compact and Hausdorff, to show that the continuous map $\pi$ is proper, it is enough to show that it is closed and that its fibers are compact.
 
 Let $C \subset \PC_{fill}(S)$ be a closed set. We need to show $\pi(C)$ is closed. So suppose $X_n \in \pi(C)$ so that $X_n \to X \in \T(S)$. We need to show that $X \in \pi(C)$. 
 
 Let $[\mu_n] \in C$ so that $\pi[\mu_n] = X_n$. Since $C \subset \PC_{fill}(S)$, these are all filling currents. Without loss of generality, their representatives in $\C_{fill}(S)$ satisfy $i(\mu_n, \mu_n) = 1$. Thus, 
 \[
  \ell_{X_n}(\mu_n) \leq 4 \sqrt{2 |\chi(S)|}
 \]
 by Lemma \ref{lem:Minimum_Length}. Lemma \ref{lem:Lipschitz_length} implies 
 \[
  i(X, \mu_n) < e^{d_{Th}(X_n, X)} i(X_n, \mu_n)
 \]
 Since $X_n \to X$, $d_{Th}(X_n, X) \to 0$. Thus, there is some $L > 0$ for which $i(X, \mu_n) \leq L$ for all $n$.
 
 The set $\{\mu \in \C(S) \ | \ i(X, \mu) \leq L \}$ is compact by \cite{Bon88}. Thus, up to passing to a subsequence, $\mu_n \to \mu$ for $\mu \in \C(S)$. 
 
 Since we assumed that $C$ is closed, we must have $\mu \in C$. In particular, $\mu$ is filling, so $\pi(\mu)$ is well defined. Moreover, for any $Y \in \T(S)$, we have 
 \[
  i(X_n, \mu_n) \leq i(Y, \mu_n)
 \]
Taking limits of both sides gives us $i(X, \mu) \leq i(Y,\mu)$. Thus, 
\[
 \pi(\mu) = X
\]
and so $X \in \pi(C)$. Therefore, $\pi(C)$ is closed.

Next, let $X \in \T(S)$. We must show that $\pi^{-1}(X)$ is compact in $\PC_{fill}(S)$. Since $\PC(S)$ is compact, it is enough to show that $\pi^{-1}(X)$ is closed in $\PC(S)$. So suppose $[\mu_n] \in \PC_{fill}$ with $\pi[\mu_n]= X$, and $[\mu_n] \to [\mu]$ in $\PC(S)$. 

By choosing appropriate representatives, we can assume that $\mu_n \to \mu$ as currents. Then we note that $X$ must be the length minimizer of $\mu$. Indeed, for all $Y \neq X \in \T(S)$, 
\[
 i(\mu_n, X) < i(\mu_n, Y)
\]
and taking limits gives us the desired inequality. By Lemma \ref{lem:No_Minimizers_Nonfilling}, non-filling currents do not have length minimizers. Thus, $\mu$ must also be filling. Therefore, $[\mu] \in \pi^{-1}(X)$, and so $\pi^{-1}(X)$ is closed inside $\PC(S)$. Thus, it is compact.  Therefore, the map $\pi$ is proper.
\end{proof}

Finally, we want to connect the systolic length of a filling current to its self-intersection number, and its length at its projection.
We will need the following remark.
\begin{rem}
 By \cite[Corollary 1.4]{BIPP}, $\Mod(S)$ acts properly discontinuously on $\C_{fill}(S)$. If we let $\Omega$ be the quotient of $\C_{fill}(S)$ by $\Mod(S)$, then $\pi$ descends to a projection 
 \[
  \pi : \Omega \to \M_g
 \]
where $\M_g$ is the moduli space of our genus $g$ surface, which is still proper.
\end{rem}
 
\begin{cor}
 For every $\epsilon$, for all $\mu\in\mathcal{C}_\mathrm{fill}(S)$ with $\pi(\mu) \in \T_\epsilon(S)$, we have
 \[
  \sqrt{i(\mu, \mu)} \asymp \sys(\mu) \asymp \ell_{\pi(\mu)}(\mu)
 \]
where we say $A \asymp B$ if there are constants $c_1, c_2>0$ so that $c_1 A \leq B \leq c_2A$, and where the constants depend only on $S$ and $\epsilon$ but not on $\mu$.
\end{cor}
\begin{proof}
	Note that the systole function descends to $\Omega = \C_{fill}(S)/\Mod(S)$, as does self-intersection. Thus, the map $f: \Omega \to  \R$ defined by 
\[
 f([\mu]) = \frac{\sys(\mu)}{\sqrt{i(\mu, \mu)}}
\]
is well-defined and continuous (Note that the self-intersection of a filling current can never be 0). 

  Next observe that both numerator and denominator in the definition of $f$ are homogeneous of degree one, and hence $f$ also descends to a continuous function on the projectivized set $\mathbb{P}\Omega = \PC_{fill}(S)/\Mod(S)$.

By the previous remark, we can also descend the map $\pi$ to a map 
 \[
  \pi : \mathbb{P}\Omega \to \mathcal M_g
 \]
which is still proper. Thus, the pre-image of the $\epsilon$-thick part of $\mathcal M_g$ is compact. Since a continuous map on a compact set achieves its maximum and minimum, there are constants $c_1= c_1(\epsilon)$ and $c_2 = c_2(\epsilon)$ for which 
\[
 c_1 \leq  \frac{\sys(\mu)}{\sqrt{i(\mu, \mu)}} \leq c_2
\]
Moreover, $c_1 > 0$ as the systolic length a filling current can never be 0.

Since the function $\mu \mapsto \ell_{\pi(\mu)}(\mu)$ is continuous and $\Mod(S)$-invariant, the same argument works for the function $g: \Omega \to \R$ with $g(\mu) = \ell_{\pi(\mu)}(\mu)/ \sqrt{i(\mu, \mu)}$, giving us the second inequality.
\end{proof}

\begin{rem}
	In fact, the same (well-known) argument from the previous proof shows that any two positive, continuous, mapping class group invariant, degree-one-homogeneous functions will be comparable on the set $\pi^{-1}\mathcal{T}_\epsilon \subset \mathcal{C}_\mathrm{fill}$. 
\end{rem}
\section{Extension of the projection}
\label{sec:Extention_to_non-filling}
By Bonahon \cite{Bon88}, there is a natural embedding of $\T(S)$ into $\PC(S)$. The set $\PC(S)$ is compact, and the closure of $\T(S)$ in this embedding is exactly the Thurston compactification. Given that the map $\pi : \PC_{fill}(S) \to \T(S)$ is proper, it is natural to ask if it can be continuously extended to a map $\PC(S) \to \overline{\T(S)}$, where $\overline{\T(S)}$ is the Thurston compactification. However, by a result of Diaz and Series, this turns out to be impossible.

Kerckhoff showed that given $\mu, \nu \in \ML(S)$ so that $\mu + \nu$ is a filling current, then $\mu + \nu$ has a unique length minimizing metric in $\T(S)$ \cite{Kerckhoff92}. Let $\mathcal A = \{\alpha_1, \dots, \alpha_n\}$ and $\mathcal B = \{\beta_1, \dots, \beta_m\}$ be two simple closed multicurves whose union fills $S$. Diaz and Series then look at the \textit{simplex of minima}, that is, the projection to $\T(S)$ of $\{\sum t_i \alpha_i + \sum s_i \beta_i \ | \ t_i, s_i > 0\}$. In \cite[Corollary 1.3]{DS03}, they show that the projection does not extend to the boundary of this set in $\PC(S)$. In fact, two different sequences that converge to $\alpha_1 + \dots + \alpha_{n-1}$ can project to sequences that converge to different elements of $\PML(S)$.

\begin{lem}[Diaz-Series, Corollary 1.3 \cite{DS03}]
 The projection $\pi: \PC_{fill}(S) \to \T(S)$ does not extend continuously to a projection from $\PC(S)$ to the Thurston compactification of $\T(S)$.
\end{lem}

On the other hand, our techniques allow us to prove an extension of Theorem 1.2 in \cite{DS03} and show that $\pi$ extends continuously at maximal, uniquely ergodic laminations. The original result is as follows. If $\mu$ and $\lambda$ were two laminations so that $\mu+ \lambda$ is a filling current, then in fact, $t \mu + (1-t)\lambda$ will be filling for all $t \in (0,1)$. Kerckhoff introduced the map $t \to \pi(t \mu + (1-t)\lambda)$. He showed that it is, in fact, a homeomorphism onto its image. This image in $\T(S)$ is called the \textit{line of minima} between $\mu$ and $\lambda$. As $t \to 0$, of course, $t \mu + (1-t)\lambda$ converges to $\lambda$. Diaz and Series show that when $\lambda$ is maximal and uniquely ergodic, then for any choice of $\mu$, the line of minima $\pi(t \mu + (1-t)\lambda)$ converges to $[\lambda] \in \PML(S)$ as $t \to 0$. We extend this result to any sequence of filling currents that converges to a maximal, uniquely ergodic lamination, and provide a converse.

\newtheorem*{thm:Uniquely_ergodic}{Proposition \ref{prop:Extension_at_Uniquely_Ergodic}}
\begin{thm:Uniquely_ergodic}
Suppose $\lambda$ is a maximal, uniquely ergodic lamination. If $\mu_n \in \C_{fill}(S)$ with $\mu_n \to \lambda$, then $\pi(\mu_n) \to [\lambda] \in \PC(S)$.

Conversely, if $\mu_n \in \C_{fill}(S)$ with $\mu_n \to \mu$ for some $\mu$, and $\pi(\mu_n) \to [\lambda]\in \PC(S)$, then, in fact, $[\mu] = [\lambda] \in \PC(S)$.
\end{thm:Uniquely_ergodic}
\begin{proof}
 Suppose we have $\mu_n \in \C_{fill}(S)$ with $\mu_n \to \mu \in \C(S)$. Let $X_n = \pi(\mu_n)$ and suppose, up to taking a subsequence, $X_n \to \lambda$ for $\lambda \in \PML(S)$. For now, we make no other assumptions on $\mu$ and $\lambda$.
 
 We can find constants $c_n$ so that $c_n X_n \to \lambda$. Since $i(\lambda, \lambda) = 0$ and $i(X_n, X_n) > 0$, we have that $c_n \to 0$. Then, by Lemma \ref{lem:Minimum_Length}, 
 \[
  i(\mu_n, c_n X_n) \leq c_n 4 \sqrt{2 |\chi(S)|} \sqrt{i(\mu_n, \mu_n)}
 \]
Taking limits, and using the fact that $c_n \to 0$ and $i(\mu_n, \mu_n) \to i(\mu, \mu)$, we get 
\[
 i(\mu, \lambda) = 0
\]
If $\mu$ is, in fact, a maximal, uniquely ergodic lamination, then the fact that $\lambda$ is a lamination implies that $[\lambda] = [\mu]$. Since every subsequence of $\pi(X_n)$ converges to $[\mu]$, the entire sequence converges as well.

On the other hand, suppose $\mu$ is some (non-filling) geodesic current, and $\lambda$ is the maximal, uniquely ergodic lamination. We'll show in this case that $[\mu] = [\lambda]$. The fact that $i(\mu, \lambda) = 0$ means that the supports of $\mu$ and $\lambda$ do not have any transverse intersection points. Thus, all the geodesics in the support of $\mu$ either lie in the support of $\lambda$, or in a complementary region of its support. But since $\lambda$ is maximal, its complementary regions are all ideal polygons. For each polygon, there are only finitely many complete geodesics contained in the interior. However, these geodesics are all isolated from the rest of the support of $\lambda$. So the only way they can be in the support of $\mu$ is if they project down to closed curves. However, since $\lambda$ is uniquely ergodic, none of its leaves are asympotic to a closed curve. Thus, the support of $\mu$ must be contained in the support of $\lambda$.  In particular, $\mu$ is also a measured lamination. Since $\lambda$ is uniquely ergodic, this means that $[\mu] = [\lambda]$ in $\PML(S)$.
\end{proof}

Although $\pi$ does not extend to a larger set inside the space of
currents, one can still wonder if there is a way to extend it (to some compactification of Teichm\"uller space). Given that $\pi$ extends at uniquely ergodic points,
and at these points the Thurston boundary agrees with the
Teichm\"uller boundary, this suggests this as a natural candidate:
\begin{ques}
	Does the map $\pi$ extend to the Teichm\"uller compactification? (or any other compactification?)
\end{ques}

 \bibliographystyle{alpha}
  \bibliography{ProjectionPropertiesArxiv}
\end{document}